\documentclass{article}

\usepackage{amssymb,amsmath,eufrak,amsthm,amscd}
\usepackage[matrix,arrow,curve]{xy}

\newtheorem{theorem}{Theorem}
\newtheorem{lemma}[theorem]{Lemma}

\newtheorem{corollary}[theorem]{Corollary}
\newtheorem{statement}[theorem]{Statement}
\newtheorem*{statement*}{Statement}
\newtheorem*{theorem*}{Theorem}
\newtheorem*{lemma*}{Lemma}
\newtheorem*{fact*}{Fact}

\newcommand{\ord}{\operatorname{ord}}

\newcommand{\Spec}{\operatorname{Spec}}
\newcommand{\QSpec}{\operatorname{QSpec}}
\newcommand{\Rad}{\operatorname{Rad}}
\newcommand{\QRad}{\operatorname{QRad}}
\newcommand{\Max}{\operatorname{Max}}
\newcommand{\QMax}{\operatorname{QMax}}

\newcommand{\Hur}{\operatorname{H}}

\theoremstyle{definition}

\newtheorem*{definition*}{Definition}
\newtheorem{example}[theorem]{Example}
\newtheorem*{example*}{Example}

\theoremstyle{remark}

\newtheorem*{remark*}{Remark}

\author{D.\,V.~Trushin}
\title{Nullstellensatz over quasifields\footnote{\Xy-pic package is used}}
\date{}

\begin{document}

\maketitle

\abstract{

We investigate the least studied class of differential rings -- the
class of differential rings of nonzero characteristic. We present
the notion of differentially closed quasifield and develop
geometrical theory of differential equations in nonzero
characteristic. The notions of quasivarieties and their morphisms
are scrutinized. Presented machinery is a basis for reduction modulo
$p$ for differential equations.

}

\section{Introduction}

We investigate the least studied class of differential rings -- the
class of differential rings of nonzero characteristic. Almost all
known technique devoted to the differential rings is based on
methods invented by Kolchin~\cite{Kl}. However, this machinery is
not appropriate in the case when the set of all differential prime
ideals is not big enough. So the rings with empty differential
spectrum are out of the field of application of the known technique.
There are a few number of early works of Yuan~\cite{Yu},
Block~\cite{Bl}, and Posner~\cite{Ps1}, \cite{Ps2}. Only at the end
of 70s William Keigher obtained first new results about simple
differential rings in nonzero characteristic. His latest
work~\cite{Ke7} is an attempt of generalization the Picard-Vessiot
theory to the case of quasifield. In our work we present a slightly
different direction: we shall generalize the algebraic geometry to
the case of quasifield. However, it should be noted that we use
Keigher's result as a background. Particulary, the notion of
quasifield  and quasispectrum are key notions in our work. The
problem of obtaining the results discussed was posed by Keigher in
April, 2007. From his point of view this work should be a basis of
the reduction modulo $p$ for differential equations. To realize our
plans we develop a new technique based on Nullstellensatz in the
case of infinite number of indeterminates and on universal property
of the Hurwitz series ring. The described method allows us not only
to define the notions we interested in but scrutinize them well.

The paper consists of the main part and appendix. In sections~2
and~3 we present the basic terms and notations. In the section~2 we
collect all terms and obvious facts. Also we give a references to
the works where discussed notions investigate more accurate. In the
section~3 we give two main statements which form a heart of our
machinery, these are theorem~\ref{GilbNull}
and~\ref{HurwitzUnivers}. In spite of the proofs of these statements
are very easy, we call them theorems because they play very impotent
role in our work. It should be noted that theorem~\ref{GilbNull} is
just useful variant of theorems proved in appendix. In the section~4
we present the notion of differentially closed quasifield and
scrutinize its structure. The main result of the section is a full
classification of differentially closed quasifields. The section is
ended by theorem~\ref{DiffClosedFull} where we collect all proved
about such quasifields facts. Section~5 is devoted to the notion of
differential closure. We define it in the same way as it done for
differential closure of the differential field of zero
characteristic~\cite[chapter~6, sec.~b]{Pz}. Also we discuss the
question of explicit construction of differential closure and the
question of its minimality theorem~\ref{ExistTheorem}
and~\ref{MinTheorem}, correspondingly. The ascending chain condition
does not hold in the ring of differential polynomials over
quasifields even for quasiradical ideals. Therefore section~6 is
devoted to its weaker analogue, the notion of $\omega$-Noetherian.
The main theorem of this section is theorem~\ref{OmegaNoeth}.
Sections from~7-th till~10-th are devoted to the geometric theory.
The notion of differentially closed quasifield allows us to define
the notion of quasivariety. Section~7 present affine quasivarieties
and its basic properties. We do not define general variant of
quasivariety but we assume that it can be done immediately by
anyone. The section mostly devoted to the translation of proved
facts to the geometric language. Additionally, we give some examples
and discuss the relation between quasivariety and spectrum. To
produce the category of quasivarieties we introduce the notion of
its morphisms. For that in section~8 we define the notion of regular
functions and mappings. Theorem~\ref{GlobSect} says that all regular
functions are polynomial (compare with differential algebraic
geometry~\cite[chapter~I, sec.~5, p.~901]{Cs}). The discussion is
finished by categorical meaning of the proved result
(theorem~\ref{Categor}). In algebraic geometry using the notion of
Noetherian, a lot of geometrical facts can be expressed in terms of
constructible sets. Noetherian does not hold for quasivarieties.
However, using more weaker notion of $\omega$-Noetherian, we are
able to prove some geometrical properties in terms of the Baire
property. Section~9 is devoted to the Baire property. We discuss the
behavior of irreducible quasivarieties under regular mappings. The
last section devoted to the notion of quasivariety is section~10, in
which the notion of dimension is introduced. The last three sections
are devoted to some technical details and generalizations. From
section~5 it follows that the structure of differential closure in
the case of quasifield with countable residue field very differs
from that in the  contrary case. Theorem~\ref{QConstr} in section~11
describe structure of differential closure in the hard case of
countable residue field. The result shows that differential closure
is constructible over basic quasifield (terms~\cite[chapter~10,
sec.~d]{Pz}). In the section~12 we generalize the obtained results
to the case of arbitrary set of differential indeterminates. The
final section is appendix. This section is mostly devoted to the
accurate proving of general form of Nullstellensatz.

\paragraph{Acknowledgement.} This text appeared due to my acquaintance with
two outstanding people and mathematicians -- William Keigher and
Yrii Razmyslov. The first one introduced me to the subject and posed
the problem, the second one provided me with the technique used
throughout the text. It gives me a great pleasure to acknowledge my
deep debt to professor Alexander V.~Mikhalev.

\section{Terms and notation}

Throughout the all text a differential ring $A$ is a commutative,
associative unitary ring with a finitely many pairwise commuting
differential operators. The set of all differential operators will
be denoted by $\Delta$ and the set of all derivations by $\Theta$.
Additionally, we assume that all ring are of nonzero characteristic.
All homomorphisms are assumed to be unitary. Differential ring will
be called simple if it has no non trivial differential ideals except
zero ideal and whole ring. All terms and notation are the same as
in~\cite{AM} and~\cite{Kl}. If something is not defined or
unexplained we mean it has the same meaning as in the mentioned
sources.

The following terms mostly are introduced in~\cite{Ke1}
and~\cite{Ke2}. We shall briefly recall them. An element $x$ of a
differential ring $A$ is called a differential nilpotent if $\theta
x$ is nilpotent whenever $\theta\in\Theta$, in other words, the
ideal $[x]$ belongs to the nilradical of the ring. An ideal $\frak
u$ will be said to be quasiradical if a quotient ring $A/\frak u$
has no nonzero differential nilpotents. The ring with the last
property will be called quasireduced. An ideal $\frak q$ will be
said to be quasiprime if it is quasiradical and primary. An ideal
$\frak m$ will be said to be quasimaximal  if it is a maximal
differential ideal.

Basic properties and relations between mentioned types of
differential ideals are scrutinized in~\cite{Ke1} and~\cite{Ke2}. We
recall just needed facts below. The sets of all radical, prime,
maximal ideals of a ring $A$ are denoted by $\Rad A$, $\Spec A$,
$\Max A$, respectively. In analogue way we shall denote by $\QRad
A$, $\QSpec A$, $\QMax A$ the set of all quasiradical, quasiprime,
and quasimaximal ideals, respectively. It is easy to see that any
quasiradical ideal can be presented as the intersection of all
quasiprime ideals containing it. The least quasiradical ideal
containing an ideal $\frak a$ will be denoted by $\frak{rad}(\frak
a)$. The last one coincides with the set of all differential
nilpotents modulo $\frak a$. In other words, the behavior of the
ideals with prefix ``quasi'' is the same as that of the ideals
without it.

Let define the following mappings on ideals:
\begin{gather*}
\pi\colon \Rad R\to \QRad R\\
\frak r\colon \QRad R\to \Rad R
\end{gather*}
where $\pi(\frak t)$ is the largest differential ideal belonging to
$\frak t$ and $\frak r(\frak u)$ is a radical of the ideal $\frak u$
in the usual sense~\cite[chapter~1, sec.~6, p.~8]{AM}. The mapping
$\pi$ is scrutinized in~\cite{Ke2}, but with the following notation
$\pi(\frak a)=\frak a_{\Delta}$. The mapping $\pi$ preserve
arbitrary intersections~\cite[prop.~1.1]{Ke2}.

The set $\Spec A$ is provided with Zarisky topology~\cite[chapter~1,
ex.~15]{AM}. Let provide topology on $\QSpec A$. The set $\QSpec A$
with topology introduced below will be called a quasispectrum. For
each subset $E\subseteq A$, let $V(E)$ denote the set of all
quasiprime ideals of $A$ which contain $E$. From
statement~\ref{QTopology} in appendix it follows that the sets
$V(E)$ satisfy the axioms for closed sets in a topological space.
For each $x\in A$, let $(\QSpec A)_x$ denote the complement of
$V(x)$ in $\QSpec A$. The mentioned sets are open and form a basis
of open sets for the introduced topology. From now and till the end
we shall mean only these two mentioned topo9logy on spectra and
quasispectra.

From the definition of the mapping $\pi$ it follows that $\pi$ is
continuous. In the case of arbitrary characteristic the mapping
$\frak r$ is not necessarily continuous, but in the case of nonzero
characteristic the following fact holds: the mappings $\pi$ and
$\frak r$ are inverse to each other bijections between the set of
all radical ideals and the set of all  quasiradical ideals.
Moreover, these two mappings provide inverse to each other
homeomorphisms between $\Spec A$ ($\Max A$) and $\QSpec A$ ($\QMax
A$) (the prove is in~\cite[lemma~13 and theorem~14]{Tr}). Let denote
everything on the following diagram
$$
\xymatrix{
    \Max R\ar@{^{(}->}[r]\ar@/^/[d]^{\pi}  &\Spec R\ar@{^{(}->}[r]\ar@/^/[d]^{\pi}    & \Rad R\ar@/^/[d]^{\pi}\\
    \QMax R\ar@{^{(}->}[r]\ar@/^/[u]^{\frak r} &\QSpec R\ar@{^{(}->}[r]\ar@/^/[u]^{\frak r}   & \QRad R\ar@/^/[u]^{\frak r}
}
$$
It should be noted, that presented topology on $\QSpec A$ coincides
with topology introduced in~\cite[sec.~2, p.~144]{Ke3}.

A ring $A$ will be said to be a quasifield  if it is simple
differential ring of nonzero characteristic. From~\cite[sec.~2,
theor.~2.1]{Yu} it follows that $A$ is a local ring of Krull
dimension zero having prime characteristic $p$. Its maximal ideal is
determined  by the following condition
$$
\frak m=\{\,x\in A\mid x^p=0\,\}.
$$
The residue field of the ideal $\frak m$ will be denoted by
$K=A/\frak m$ and called a residue field of quasifield $A$. An
additional information about quasifield structure can be found
in~\cite[sec.~2, prop.~2.8]{Ke2}.

\section{Basic technique}

The most popular technique for investigating differential rings is
based on characteristic sets. However, this machinery does not
provide desired results in nonzero characteristic. The section is
devoted to the adequate ring machinery providing our needs.

The following result was found in~\cite[sec.~2, ex.~2.1, p.~18]{Yr}.
However, its original prove is based on algebras representation
theory. We give the most useful version of the result. The prove of
the full sequence of the general Hilbert Nullstellensatz is given in
appendix.

\begin{theorem}\label{GilbNull}
Let $K\subseteq L$ be a field extension such that $L$ is generated
over $K$ by not more then $\kappa$ elements ($\kappa$ is an
arbitrary cardinal). Then, if $|K|>\kappa$ then $L$ is algebraic
over $K$.
\end{theorem}
\begin{proof}
The proof follows from~\ref{WeakNullAlg}.
\end{proof}

From the result above the following one immediately follows.

\begin{corollary}\label{NullCoroll}
Let $K$ be an algebraically  closed field and $X$ be arbitrary set
such that $|X|<|K|$. Then for each maximal ideal $\frak m$ of the
ring $K[X]$ its quotient ring $K[X]/\frak m$ coincides with $K$.
\end{corollary}

Consider the example showing that the condition $|X|<|K|$ can not be
omitted. This example is a partial case of
theorem~\ref{WeakNullAlg}.

\begin{example}\label{NullExample}
Let $K$ be a countable algebraically closed field and let $L$ be
algebraical closure of $K(x)$, where $x$ is transcendental over $K$.
It is easy to see that $L$ is countable, and therefore there is a
surjective homomorphism of the ring $K[x_n]_{n\in\mathbb N}$ onto
$L$. The kernel of the last one is the desired example.
\end{example}

For any commutative ring $A$ the ring of Hurwitz series is defined
in~\cite[sec.~2]{Ke5}. We shall define it by $\Hur A$. The Hurwitz
series ring is a differential ring. From the described
correspondence between spectrum and quasispectrum it follows that
$\Hur A$ is a quasifield iff $A$ is a field. The set of all series
with zero free term will be denoted by $\Hur A_1$. In nonzero
characteristic the ideal $\Hur A_1$ consists of nilpotent elements.
The mapping from $\Hur A$ to $A$, maps each series to its free term,
is a ring homomorphism and is denoted by $\pi$. The following result
is proved in~\cite[sec.~2, prop.~2.1]{Ke5} and~\cite[p.~100]{Ke6}
(in the case of one derivation) and in~\cite[p. 216]{Tr} (in
generale case).

\begin{theorem}\label{HurwitzUnivers}
Let $A$ be an arbitrary ring, $B$ be a differential ring, and
$\varphi\colon B\to A$ be a ring homomorphism. Then there exists a
unique differential ring homomorphism $\Phi$ (the Taylor
homomorphism) such that the following diagram is commutative
$$
\xymatrix{
                                       & \Hur A\ar[d]^{\pi}  \\
 B\ar[r]^{\varphi}\ar@{-->}[ur]^{\Phi} & A
}
$$
\end{theorem}

For any quasifield $Q$, its residue field $K$, and arbitrary set $Y$
of differential indeterminates  over $Q$ the quotient ring of
$Q\{Y\}$ by its nilradical will be denoted by $K\{Y\}$.

\section{Differentially closed quasifields}

We shall introduce the notion of differentially close quasifield and
classify all such quasifields up to isomorphism.

Let $Q$ be any quasifield. Consider the ring of differential
polynomials $Q\{y_1,\ldots,y_n\}$ over $Q$. For each subset
$E\subseteq Q\{y_1,\ldots,y_n\}$, let $V(E)$ denote the set of all
common zeros in $Q^n$ for the polynomials in $E$. Conversely, for
each subset $X$ of $Q^n$, let $I(X)$ denote the set of all
polynomials vanishing on $X$. It is clear, that for any differential
ideal $\frak a$ there is the inclusion $\frak{rad}(\frak a)\subseteq
I(V(\frak a))$.

A quasifield $Q$ will be said to be differentially closed if for any
natural number $n$ and any differential ideal $\frak a$ of
$Q\{y_1,\ldots,y_n\}$ there is the equality $\frak{rad}(\frak
a)=I(V(\frak a))$.

\begin{statement}\label{BasicDef}
The following conditions on quasifield $Q$ are equivalent:
\begin{enumerate}
\item $Q$ is differentially closed.
\item For any natural number $n$ and any proper differential ideal
$\frak a$ of  $Q\{y_1,\ldots,y_n\}$ the set $V(\frak a)$ is not
empty.
\item Any differentially finitely generated over $Q$ quasifield coincides with $Q$.
\end{enumerate}
\end{statement}
\begin{proof}
Implication (1)$\Rightarrow$(2) is obvious.

(2)$\Rightarrow$(3). Let $L$ be differentially finitely generated
over $Q$ quasifield. Then $L$ is of the form
$Q\{y_1,\ldots,y_n\}/\frak m$ for some quasimaximal ideal $\frak m$.
By the data, the set $V(\frak m)$ contains a point $a\in Q^n$. Then
$\frak m\subseteq I(a)$ and hence $\frak m=I(a)$, thus $L=Q$.

(3)$\Rightarrow$(1). Let $\frak a$ be a differential ideal of
$Q\{y_1,\ldots,y_n\}$ and let $f$ is not in $\frak r(\frak a)$. Then
there exists a $\theta\in \Theta$ such that $\theta f$ is not
nilpotent modulo $\frak a$. Therefore, algebra
$B=Q\{y_1,\ldots,y_n\}_{\theta f}/\frak a$ is not zero. Consider a
quasimaximal ideal $\frak m$ in $B$, then $B/\frak m=Q$. The images
of the elements $y_1,\ldots,y_n$ give us a point $a$ in $Q^n$
belonging to $V(\frak a)$ and such that $f(a)\neq 0$.
\end{proof}

The following statements guaranties that differentially closed
quasifields exist.

\begin{statement}\label{ClosedHurwitz}
Let $K$ be more than countable algebraically closed field. Then
$\Hur K$ is a differentially closed quasifield.
\end{statement}
\begin{proof}

Let show that the condition~(3) of statement~\ref{BasicDef} is
satisfied. Indeed, let $L$ be a differentially finitely generated
over $\Hur K$ quasifield. Then it can be presented in the following
form
$$
L=\Hur K\{y_1,\ldots,y_n\}/\frak m,
$$
where $\frak m$ is a quasimaximal ideal. We shall construct a
homomorphism from $L$ to $K$. Let $\frak n$ be the radical of the
ideal $\frak m$. It is easy to see, that it contains the ideal $\Hur
K_1\{y_1,\ldots,y_n\}$. Then from corollary~\ref{NullCoroll} it
follows that
$$
\Hur  K\{y_1,\ldots,y_n\}/\frak n= K\{y_1,\ldots,y_n\}/\frak n= K.
$$
Let denote the constructed homomorphism by $\varphi$. Applying
theorem~\ref{HurwitzUnivers} we get that there exists a unique
homomorphism $\Phi$ such that
$$
\xymatrix{
                                       & \Hur K\ar[d]^{\pi}  \\
L\ar@{-->}[ur]^{\Phi} \ar[r]^-{\varphi} & K }
$$
From the uniqueness of the Taylor homomorphism it follows that
$\Phi$ is a homomorphism over $\Hur K$. Since $L$ is quasifield, the
mapping $\Phi$ is an isomorphism over $\Hur K$.
\end{proof}

Example~\ref{NullExample} implies the following result.

\begin{statement}\label{ResidueClosed}
Let $Q$ be differentially closed quasifield and let $K$ be its
residue field. Then $K$ is more than countable.
\end{statement}
\begin{proof}
Indeed, suppose that contrary holds. Consider the ring $K[\theta
y]_{\theta\in \Theta}$ and its ideal $\frak m$ constructed as in
example~\ref{NullExample}. The ideal can be extended to the ideal
$\frak m'$ of $Q\{y\}$. Let put $\frak q=\pi(\frak m')$. We shall
shaw that quasifield $D=Q\{y\}/\frak q$ does not coincide with $Q$.
Let compare residue fields for that.
$$
Q\{y\}/\frak m'=K[\theta y]_{\theta\in \Theta}/\frak
m=\overline{K(x)}.
$$
The contradiction finished the proof.
\end{proof}

\begin{statement}\label{ClassClosedHurwitz}
Quasifield $\Hur K$ is differentially closed iff $K$ is more than
countable algebraically closed field.
\end{statement}
\begin{proof}
The right to left statement is proved in~\ref{ResidueClosed}. The
other side statement follows from the last one.
\end{proof}

\begin{statement}\label{DiffClosSingle}
Let $Q$ be a quasifield with a residue field $K$. Suppose that any
quasifield $L$, generated by one element over $Q$, coincides with
$Q$. Then the Taylor homomorphism is an isomorphism between $Q$ and
$\Hur K$.
\end{statement}
\begin{proof}
Let $\varphi\colon Q\to K$ be a quotient homomorphism and let
$\Phi\colon Q\to\Hur K$  be corresponding to $\varphi$ the Taylor
homomorphism. Since $Q$ is a quasifield, we only need to check that
$\Phi$ is surjective. Let suppose that contrary holds. In other
words, there exists an element $\eta\in\Hur K$ not belonging to $Q$.
Let denote its coefficients by $a_\theta$ and consider the ideal
$(\ldots,\theta y-a_\theta,\ldots)$ in $K\{y\}$. This ideal extends
to maximal ideal $\frak m$ in $Q\{y\}$. Let $\frak q$ be
corresponding to $\frak m$ quasimaximal ideal. Then, from the one
hand, $Q\{y\}/\frak q$ coincides with $Q$ by the data. From the
other hand, the image of $Q\{y\}/\frak q$ under the Taylor
homomorphism contains $\eta$, contradiction.
\end{proof}

So, we have just classified all differentially closed quasifields.
Let underline it in the following statement.

\begin{statement}\label{DiffClosedClass}
Differentially closed quasifields are exactly the Hurwitz series
rings over more than countable algebraically closed fields.
\end{statement}

Now we are able to prove generale result about differentially closed
quasifields.

\begin{theorem}\label{DiffClosedFull}
The following conditions on quasifield $Q$ are equivalent:
\begin{enumerate}
\item $Q$ is differentially closed.
\item For any natural number $n$ and any proper differential ideal $\frak a$ of $Q\{y_1,\ldots,y_n\}$
the set $V(\frak a)$ is not empty.
\item Any differentially finitely generated over $Q$ quasifield coincides with $Q$.
\item Any quasifield generated over $Q$ by one single element
coincides with $Q$.
\item Any quasifield generated over $Q$ by not more than countable set of elements
coincides with $Q$.
\item $Q$ isomorphic to the Hurwitz series rings over more than countable algebraically closed fields.
\end{enumerate}
\end{theorem}
\begin{proof}
Only one thing we need to show is that condition~(5) follows from
the others. From theorem~\ref{HurwitzUnivers} it follows that we
need to show that for any such quasifield $L$ there is a
homomorphism from $L$ to the residue field of $Q$. Let the residue
field of $Q$ be denoted by $K$. The residue field of $L$ is not more
than countably generated over $K$ and thus from
corollary~\ref{NullCoroll} it coincides with $K$, q.e.d.
\end{proof}

\section{Differential closure}

The following step is to define a differential closure for arbitrary
quasifield. Let $Q$ be a quasifield. We shall say that $L$ a
differential closure of $Q$ if for any differentially closed
quasifield $D$ containing $Q$ there is a homomorphism from $L$ to
$D$ over $Q$. We shall prove that differential closure exists and
unique up to isomorphism. Additionally we shall discuss the
minimality question.

Let us describe the construction of differential closure. Let $Q$ be
a quasifield with a residue field $K$. If the field $K$ is more than
countable, then $L$ will be defined as an algebraic closure of $K$.
If the field $K$ is not more than countable, then  $L$ will be
defined as an algebraic closure of
$K(x_\alpha)_{\alpha\in\omega_1}$, where $\omega_1$ is the first
more than countable cardinal. The differential closure
$\overline{Q}$ will be defined as $\Hur L$. Moreover, we have the
following sequence of homomorphisms
$$
Q\to K\to L.
$$
Applying the Taylor homomorphism to the last one, we get the
inclusion of $Q$ to its differential closure $\overline{Q}$.

\begin{theorem}\label{ExistTheorem}
Let $Q$ be a quasifield with a residue field $K$. Then a
differential closure of $Q$ is unique up to isomorphism and is
isomorphic to the quasifield $\overline{Q}$.
\end{theorem}
\begin{proof}
The fact that $\overline{Q}$ is differentially closed follows from
statement~\ref{ClassClosedHurwitz}. We shall prove the universal
property. Let $D$ be a differentially closed quasifield containing
$Q$ and let the residue field of $D$ is denoted by $F$. Then $K$ is
embedded to $F$. Moreover, from statement~\ref{ClassClosedHurwitz}
it follows that $F$ is more than countable algebraically closed
field. Therefor, the embedding $K$ to $F$ can be extended to
embedding $L$ to $F$ ($L$ is a residue field of $\overline{Q}$).
Theorem~\ref{HurwitzUnivers} guaranties that there is an embedding
$\overline{Q}$ to $D$. From the uniqueness of the Taylor
homomorphism it follows that this embedding is over $Q$.
\end{proof}

\begin{statement}\label{SubQuasiState}
Let $Q$ be a quasifield with a residue field $K$. If the field $K$
is more than countable then for any differential ring $D$ such that
$Q\subseteq Q\subseteq \overline{Q}$ it follows that $D$ is a
quasifield.
\end{statement}
\begin{proof}
To prove the result we just need to show that for each element
$\eta\in \overline{Q}$ the ring $Q\{\eta\}$ is a quasifield. From
the definition we have $\overline{Q}=\Hur \overline{K}$. Consider
the element $\eta\in \Hur \overline{K}$ and let $a_\theta$  be its
coefficients. Let define the ideal $(\ldots,\theta
y-a_\theta,\ldots)$ in $\overline{K}\{y\}$. Since $\overline{K}$ is
integral over $K$, then the mentioned ideal contract to a maximal
ideal $\frak m'$ in $K\{y\}$. The ideal $\frak m'$  extends to the
maximal ideal $\frak m$ in $Q\{y\}$. Let $\frak q=\pi(\frak m)$.
Then, from one hand, $Q\{y\}/\frak q$  is a quasifield. From the
other hand, the ring $Q\{\eta\}$ coincides with the image of
$Q\{y\}/\frak q$ under the Taylor homomorphism.
\end{proof}

\begin{theorem}\label{MinTheorem}
Let $Q$ be a quasifield with a residue field $K$. Differential
closure $\overline{Q}$ is minimal over $Q$ iff the field $K$ is more
than countable.
\end{theorem}
\begin{proof}
If the field $K$ is more than countable then the desired result
follows from the previous statement. Indeed, suppose that contrary
holds, that there is a differentially closed quasifield $D$ with
condition $Q\subseteq D\subseteq \overline{Q}$. But for each element
$y\in \overline{Q}\setminus D$ the algebra $D\{y\}$ is a quasifield
and does not coincide with $D$, contradiction.

Let $K$ be not more than countable. From the definition of the
quasifield $\overline{Q}$ it follows that there is a subfield $L'$
in $L$ containing $K$ and isomorphic to $L$. Thus the quasifield
$\Hur L'$ is a differentially closed quasifield such that
$Q\subsetneq\Hur L'\subsetneq \overline{Q}$.
\end{proof}

\section{$\omega$-Noetherian}

We shall say that a particulary ordered set $S$ is
$\omega$-noetherian ($\omega$ is a countable cardinal) if each more
than countable ascending chain is stable, or, in other words, every
strictly ascending chain of elements is not more than countable. A
ring $A$ will be said to be $\omega$-noetherian if the set of all
its ideals is $\omega$-noetherian. The following statement is
obvious.

\begin{statement}
A ring $A$ is $\omega$-noetherian iff each ideal of $A$ is not more
than countably generated.
\end{statement}

The main example of an $\omega$-noetherian ring is a countably
generated algebra over a field. It follows from the next result.

\begin{statement}\label{OmegaNoetherianSt}
The ring of polynomials $\mathbb K[x_n]_{n\in\mathbb N}$ is
$\omega$-noetherian.
\end{statement}
\begin{proof}

For any ideal $\frak a$ we have $\frak a=\cup_n\frak a_n$, where
$$
\frak a_n=\frak a\cap K[x_1,\ldots,x_n]
$$
But each ideal $\frak a_n$ is finitely generated. Therefor, $\frak
a$ is not more than countably generated.
\end{proof}

The following fact is the analogue of Ritt-Raudenbusch's basis
theorem.

\begin{theorem}\label{OmegaNoeth}
Let $Q$ be a quasifield with a residue field $K$. Then the set
$\QRad Q\{y_1,\ldots,y_n\}$ is $\omega$-noetherian.
\end{theorem}
\begin{proof}
Let $\frak n$ be the nilradical of $Q$. Then the particulary ordered
sets

$$
\QRad Q\{y_1,\ldots,y_n\}\:\mbox{ and }\:\Rad Q\{y_1,\ldots,y_n\}
$$
are isomorphic. The last one is isomorphic to
$$
\Rad Q/\frak n\{y_1,\ldots,y_n\}=\Rad K\{y_1,\ldots,y_n\}.
$$
Statement~\ref{OmegaNoetherianSt} finished the proof.
\end{proof}

\begin{corollary}\label{GeneratorCor}
Each quasiradical ideal $\frak u$ of $Q\{y_1,\ldots,y_n\}$ is not
more than countably generated as a quasiradical ideal.
\end{corollary}

\section{Quasivarieties}

We shall fix a differentially closed quasifield $Q$ until the end of
the text. Its residue field will be denoted by $K$. We shall deal
with differential polynomial ring $Q\{y_1,\ldots,y_n\}$ and its
quotient by nilradical ring $K\{y_1,\ldots,y_n\}$.

Statement~\ref{DiffClosedClass} implies that
$$
Q^n=(\Hur K)^n=K^\Theta\times\ldots \times K^\Theta.
$$
Since $\QSpec Q\{y_1,\ldots, y_n\}$ is homeomorphic to $\Spec
K\{y_1,\ldots,y_n\}$ we are able to identify the following
topological spaces
$$
Q^n=\QMax Q\{y_1,\ldots, y_n\}=\Max
K\{y_1,\ldots,y_n\}=K^\Theta\times\ldots \times K^\Theta,
$$
where the element $(a_{\theta_1,\,1},\ldots,a_{\theta_n,\,n})$ in
$K^\Theta\times\ldots \times K^\Theta$ corresponds to the element
$(q_1,\ldots,q_n)$  in $Q^n$ such that $a_{\theta_i,\,i}$ are the
coefficients of the series $q_i$.

Consider an affine space $Q^n$. For each set of differential
polynomials $E$ in $Q\{y_1,\ldots,y_n\}$, let define  the following
subset
$$
V(E)=\{\,x\in Q^n\mid \forall g\in E:\:g(x)=0\,\}.
$$
A subset $X$ in $Q^n$ will be said to be an affine quasivariety (or
simply quasivariety) if $X$ is of the form $V(E)$ for some $E$. The
sets $V(E)$ satisfy the axioms for closed sets in a topological
space. The mentioned topology coincides with topology on spectrum
under the bijection above.

Let $X$ be an arbitrary subset in $Q^n$. Then we shall define the
ideal
$$
I(X)=\{\,f\in Q\{y_1,\ldots,y_n\}\mid f|_X=0\,\}.
$$
From the definition of differentially closed quasifield the
following result follows immediately.

\begin{statement}
The mappings $I$ and $V$ are inverse to each other bijections
between the set of all quasivarieties in $Q^n$ and the set of all
quasiradical ideals in $Q\{y_1,\ldots,y_n\}$.
\end{statement}

It should be noted that in contrast to algebraic geometry
(differential algebraic geometry) the set of all irreducible
components of quasivariety is not necessarily finite.

\begin{example}
The notion of irreducible component is a topological notion. Since
every quasivariety is homeomorphic to a spectrum of an algebra that
is countably generated over a field, we just need to produce an
example of such an algebra. Consider the polynomial ring
$K[x_n]_{n\in\mathbb N}$ and the ideal
$$
\frak a=(\ldots,x_i x_j,\ldots)_{i\neq j}.
$$
It is easy to see that the ideals
$$
\frak p=(\ldots,x_j,\ldots)_{j\neq l}
$$
are its minimal prime components.
\end{example}

We shall call a topology space $\omega$-noetherian if its every more
than countable descending chain of closed sets is stable. In other
words, statement~\ref{OmegaNoeth} means the following.

\begin{statement}
Every quasivariety is $\omega$-noetherian.
\end{statement}

Another variant of the same result is the following.

\begin{statement}
Every system of differential equations over quasifield is equivalent
to its not more than countable subsystem.
\end{statement}

Any element of $Q\{y_1,\ldots,y_n\}$ determines a differential
polynomial function $Q^n\to Q$. Let $X$ be a quasivariety in $Q^n$.
A restriction of a differential polynomial mapping to $X$ is called
a differential polynomial function on $X$. The set of all
differential polynomial functions on $X$ will be denoted by
$Q\{X\}$. This ring is isomorphic to $Q\{y_1,\ldots,y_n\}/I(X)$. The
quotient ring of $Q\{X\}$ by its nilradical will be denoted by
$K\{X\}$. We know that $X$ coincides with $\QMax Q\{X\}$ and $\Max
K\{X\}$. A principal open set is a set of the following form
$$
X_g=\{\,a\in X\mid g(a)\neq 0\,\}.
$$
It is easy to see that homeomorphism between $X$ and $\QMax Q\{X\}$
preserves the notion of principal open. But it should be noted that
this notion is not preserved under homeomorphism between $\QMax
Q\{X\}$ and $\Max K\{X\}$. Namely, principal open set $(\Max
K\{X\})_g$ is mapped to principal open of the form $(\QMax
Q\{X\})_{g^p}$. The principal open sets of the form $(\QMax
Q\{X\})_{g^p}$ will be called a ``good'' principal open sets. As we
can see, there is the equality
$$
X_{g^p}=\QMax(Q\{X\}_g).
$$

\begin{statement}
Let $A$ be a ring of nonzero characteristic $n$. Then ``good'' principal
open sets form a basis of topology on $\QSpec A$ ($\QMax A$).
\end{statement}
\begin{proof}
We shall prove it for quasispectrum. Let the following notation be
fixed $X=\Spec A$ and $Y=\QSpec A$. Let $\pi\colon X\to Y$ be the
mentioned homeomorphism. Consider
\begin{gather*}
Y_f=Y\setminus V([f])=\pi(X\setminus V(\sum_k
(f^{(k)})))=\pi(X\setminus \bigcap_k V(f^{(k)}))=\\
=\pi(\bigcup_k X_{f^{(k)}})=\bigcup_k Y_{(f^{(k))^n}}
\end{gather*}
Q.E.D.
\end{proof}

The geometric version of previous lemma is

\begin{corollary}\label{GoodTopologyBasis}
The ``good'' principal open sets form a basis of topology on $X$.
\end{corollary}

\section{Regular functions and a structure sheaf}

Let $X$ be a quasivariety and let $f\colon X\to Q$ be a function on
$X$. We shall say (compare with~\cite[sec.~3]{Ke3}) that $f$ is
regular in $x\in X$ if there are an open neighborhood $U$ containing
$x$ and elements $h,g\in R$ such that for every element $y\in U$
$g(y)$ is not nilpotent and $f(y)=h(y)/g(y)$. The condition on $g$
can be stated as follows: for each element $y\in U$ $g^p(y)\neq0$.
For any subset $Y$ if $X$ a function is said to be regular on $Y$ if
it is regular in each point of $Y$. The set of all regular functions
on open subset $U$ of $X$ will be denoted by $\mathcal O_X(U)$.
Since definition of regular function arises from local condition the
set of rings $\mathcal O_X$ form a sheaf that will be called a
structure sheaf on $X$. This definition coincides with the usual one
in algebraic geometry. From the definition it follows that there is
the inclusion $Q\{X\}\subseteq \mathcal O_X(X)$. The very important
fact that the other one is also true.

\begin{theorem}\label{GlobSect}
For arbitrary quasivariety $X$ there is the equality
$$
Q\{X\}=\mathcal O_X(X).
$$
\end{theorem}
\begin{proof}
Let $f$ be a regular function on $X$. From the definition of regular
function it follows that for each point $x\in X$ there exist a
neighborhood $U_x$ and elements $h_x,g_x\in Q\{X\}$ such that for
every element $y\in U_x$ $g^p_x(y)\neq 0$ and $f(y)=h_x(y)/g_x(y)$.
From corollary~\ref{GoodTopologyBasis} it follows that each $U_x$
can be presented as a union of ``good'' principal open sets and this
family cover $X$. Since $X$ is homeomorphic to a prime spectrum it
follows that $X$ is compact and thus there is a finite family of
elements $\{\,a_i,b_i,c_i\,\}$ in $Q\{X\}$ such that $X_{c^p_i}$
cover $X$ and for each $y\in X_{c^p_i}$ $b^p_i(y)\neq0$ and
$f(y)=a_i(y)/b_i(y)$. Since $b_i$ is not nilpotent in $X_{c^p_i}$ we
have $X_{(c_ib_i)^p}$ coincides with $X_{c^p_i}$. Thus we can
suppose that there is a finite family of elements $\{\,a_i,b_i\,\}$
in $Q\{X\}$ such that  $X_{b^p_i}$ cover $X$ and for each $y\in
X_{b^p_i}$ $b^p_i(y)\neq0$ and
$$
f(y)=\frac{a_i(y)}{b_i(y)}=\frac{a_i(y)b^{2p-1}_i(y)}{b^{2p}_i(y)}=\frac{a'_i(y)}{b'_i(y)}.
$$
where $b'_i$ are constant and $a'_i(y)=0$ whenever $b'_i(y)=0$.
Since $X$ is covered by $X_{b'_i}$ we have
$[b'_1,\ldots,b'_n]=(b'_1,\ldots,b'_n)=(1)$. Therefore
$$
1=h_1b'_1+\cdots+h_n b'_n
$$
And the element $r=h_1a'_1+\cdots+h_n a'_n$ coincides with $f$
everywhere on $X$.
\end{proof}

The following result says that the structure sheaf is really a
structure sheaf.

\begin{statement}
For any open subset $U\subseteq X$ the mapping
$$
\QMax \mathcal O_X(U)\to \QMax Q\{X\}
$$
is a homeomorphism onto its image $U$. In other words, there is a natural homeomorphism
$$
\QMax \mathcal O_X(U)=U.
$$
\end{statement}
\begin{proof}
The set $U$ can be presented as a union of ``good'' open subsets
$X_{g^p_i}$. Let denote $\QMax \mathcal O_X(U)$ by $Y$. As we can
see $Y=\cup_i Y_{g^p_i}$. It is suffices to show that for each $i$
and $j$ the mentioned mapping induces a homeomorphism between
$X_{g^p_i}\cap X_{g^p_j}$ and $Y_{g^p_i}\cap Y_{g^p_j}$.

Indeed, consider a sequence of homomorphisms
$$
Q\{X\}\to\mathcal O_X(U)\to Q\{X_{(g_ig_j)^p}\}.
$$
Localizing by $(g_i g_j)$ we have
$$
Q\{X\}_{g_ig_j}\to\mathcal O_X(U)_{g_ig_j}\to Q\{X\}_{g_ig_j}.
$$
Everything that we need to check is that the second arrow in the
last diagram is injective. But its kernel is a localization of that
of $\mathcal O_X(U)\to Q\{X_{(g_ig_j)^p}\}$. The last one consists
of all functions $f$ that are zero on $X_{g^p_i}\cap X_{g^p_j}$ and
thus satisfy the condition $(g_ig_j)^p f =0$.
\end{proof}

Let denote $\mathcal O_X(U)$ by $Q\{U\}$. Only fact needed to
remember is that $Q\{U\}$ is not necessary differentially finitely
generated over $Q$ for arbitrary open $U$. And it should be noted
that
$$
Q\{X_f\}\neq Q\{X\}_f
$$
for a nonconstant ($\delta f=0$ for all $\delta\in\Delta$) element $f$ (for example
for any nonconstant nilpotent in $Q$).

Form theorem~\ref{GlobSect} it follows that a function coinciding
with a restriction of some differential polynomial on $X$ and a
regular function are the same things. Now we shall define a
morphisms. Let $X$ and $Y$ be quasivarieties such that $X\subseteq
Q^n$ and $Y\subseteq Q^m$. The mapping $\varphi\colon X\to Y$ will
be said to be a regular mapping if it can be defined by polynomials
$\varphi_k(x_1,\ldots,x_n)$ as follows
$$
\varphi=(\varphi_1(x_1,\ldots,x_n),\ldots,\varphi_m(x_1,\ldots,x_n))
$$
It is clear that a regular mapping $\varphi\colon X\to Y$ defines a
differential homomorphism
$$
\varphi^*\colon Q\{Y\}\to Q\{X\}
$$ by the formula
$\xi\mapsto \xi\circ \varphi$. Conversely, any differential
homomorphism $\phi\colon Q\{Y\}\to Q\{X\}$ define a regular mapping
$$
\phi^*\colon \QMax Q\{X\}\to \QMax Q\{Y\}
$$
by the rule $\phi^*(\frak q)=\phi^{-1}(\frak q)$, i.e. $\phi^*\colon X\to Y$.
The last one implies the following result.

\begin{theorem}\label{Categor}
There is a contravariant equivalence between the category of
quasireduced differentially finitely generated algebras over $Q$ and
the category of quasivarieties over $Q$.
\end{theorem}

\section{The Baire property}

In algebraic geometry properties of the morphisms can be explained in terms of constructible sets.
The section is devoted to an adequate analogue of this notion. We shall demonstrate the technique
analogue to the Noether's normalization and derive from that different geometrical applications.

We shall say that the Baire property holds for a topological space
$X$ if the intersection of any countable family of dense open
subsets is not empty. A subset is said to be big enough if it
contains a mentioned intersection.

\begin{lemma}
A ``good'' principal open set $X_{f^p}$ is dense in $X$ iff $f$ is
not a zero-divisor in $K\{X\}$.
\end{lemma}
\begin{proof}
The proof immediately follows from the fact that $\Max A_f$ is dense
iff $f$ is not zero-devisor in $A$.
\end{proof}

\begin{theorem}
The Baire property holds for irreducible quasivarieties.
\end{theorem}
\begin{proof}

In irreducible variety any non empty open subset is dense. Since ``good'' principal
open subsets form a basis of topology, we only need to prove the fact for
``good'' principal open subsets.

Let $X_{(s_n)^p}$ be a family of ``good'' principal open sets. From the previous
statement it follows that $s_n$ is not a zero devisor in $K\{X\}$. Let $S$ be a
multiplicative set generated by all elements $s_n$. Then algebra $S^{-1}Q\{X\}$ be
nonzero and differentially not more then countably generated over $Q$. Therefor, from
theorem~\ref{DiffClosedFull} item~(5) it follows that for any quasimaximal idial $\frak m$
in $S^{-1}Q\{X\}$ there is the equality
$$
S^{-1}Q\{X\}/\frak m=Q.
$$
Thus $\QMax S^{-1}Q\{X\}$ corresponds to the intersection  $\cap_n X_{(s_n)^p}$ and is not empty.
\end{proof}

The following result will be proved using the same ideas.

\begin{theorem}
Let $\varphi\colon X\to Y$ be a regular mapping with a dense image
and $Y$ is irreducible. Then the image $\varphi$ is a big enough
subset in $Y$.
\end{theorem}
\begin{proof}

Since image of $X$ is dense in $Y$ the algebra $Q\{Y\}$ is embedded in $Q\{X\}$.
From the correspondence between maximal spectrum and quasimaximal spectrum it
follows that we need to prove the statement for the pair $K\{Y\}\subseteq K\{X\}$.
The algebra $K\{X\}$ is not more than countable union of finitely generated algebras
over $K\{Y\}$, in other words
$$
K\{X\}=\bigcup_{n=0}^\infty A_k,
$$
where $A_k$ is finitely generated over $K\{Y\}$. For each $k$ there
is $s_k$ in $K\{Y\}$ such that the mapping $\Max (A_k)_{s_k}\to \Max
K\{Y\}_{s_k}$ is surjective (\cite[chapter~5, ex.~21]{AM} is a some
variation on the Neother's normalization theme). Let $S$ be a
multiplicative set generated by all elements $s_k$. Then the mapping
$$
\Max S^{-1}K\{X\}\to \Max S^{-1} K\{Y\}
$$
is surjective and the set $\Max S^{-1} K\{Y\}$ corresponds to $\cap_{k}Y_{s_k^p}$ and belongs to image of $\varphi$.
\end{proof}

\section{Dimension}

There is a standard technique of measuring a filtered algebra.
We shall apply this machinery to the case of quasivariety.

Let $X$ be a quasivariety. Consider the ring $W=K\{X\}$. The ring
$W$ is generated by the elements $\theta y_k$. We can define the
subalgebra
$$
W_r=K[\theta_1 y_1,\ldots,\theta_n y_n\mid \ord \theta_k\leqslant r]
$$
and the function $\omega_X(t)=\dim W_t$ (Krull dimension). Additionally, we consider the
following equivalence relation
$$
f(n)\sim g(n)\,\,\Leftrightarrow \,\, f(n)\leqslant g(n+n_0) \mbox{
è } g(m)\leqslant f(m+m_0)
$$

\begin{statement}
The equivalence class of the function $\omega_X$ does not depend on the choice of differential generators.
\end{statement}
\begin{proof}
Let $W_r$ and $W'_r$ be two different filtrations on $W$ generated
by two families of generators. Since $W$ is covered by both
sequences of subalgebras $W_r$ and $W'_r$ there are two numbers
$m_0$ and $n_0$ such that   $W_r\subseteq W'_{r+m_0}$ and
$W'_r\subseteq W_{r+n_0}$.
\end{proof}

The mentioned equivalence class we shall call a dimension of quasivariety.
It should be noted that $\omega_X(t)$ may not coincide with polynomial
function even for sufficiently large $t$.

\section{Well-ordering on differential closure}

Now let $Q$ be an arbitrary quasifield with a residue field $K$. It
is clear that statement~\ref{SubQuasiState} does not hold in the case
of not more than countable residue field. For example, the residue
field of $\overline{Q}$ is not algebraic over $K$, and thus there is a
subring in the residue field of $\overline{Q}$ that is not a field. We
shall shaw that even in this case there is some weaker analogue of
statement~\ref{SubQuasiState} holds.

\begin{statement}
Let $L_1$ and $L_2$ be quasifields containing $Q$. Then
$$
L_1\otimes_Q  L_2\neq0.
$$
\end{statement}
\begin{proof}
Let $K_1$ and $K_2$ be the residue fields of $L_1$ and $L_2$ respectively.
Then the desired result immediately follows from the following
$$
L_1\otimes_Q  L_2\to K_1\otimes_{K} K_2\neq0.
$$
\end{proof}

\begin{theorem}\label{QConstr}
Let $Q$ be a quasifield with a residue field $K$. Then there is a
well-ordered chain $\{\,L_\alpha\,\}$ of quasisubfields of
$\overline{Q}$ such that $\emph{(i)}$ $L_0=Q$, $\emph{(ii)}$
$L_{\alpha+1}$ is generated by one single element over $L_\alpha$,
$\emph{(iii)}$ $L_\beta=\cup_{\alpha<\beta} L_\alpha$ for limit
ordinals.
\end{theorem}
\begin{proof}

We shall organize the proof in the following manner. Firstly, we shall
construct a quasifield with desired properties and, secondary, we shall shaw,
that the constructed quasifield coincides with $\overline{Q}$.

Let $\{\,B_\alpha\,\}_{\alpha\in \Lambda}$ be the family of all differentially
finitely generated quasifields over $Q$ up to isomorphism. We may suppose that
$\Lambda$ is well-ordered. Using induction by $\alpha$ we shall construct a
sequence of quasifields $L^1_\alpha$. The quasifield $L^1_0$ will be defined as $Q$.
For any limit ordinal  $\alpha$ we have
$$
L^1_\alpha=\cup_{\beta<\alpha}L^1_\beta
$$
and
$$
L^1_{\alpha+1}=(L^1_{\alpha}\otimes_Q B_{\alpha+1})/\frak m,
$$
for some quasimaximal ideal $\frak m$. At the last step we shall
construct a quasifield $L^1$. This quasifield  satisfies the
following two properties: any homomorphism of $Q$ to any differentially
closed quasifield can be extended to a homomorphism of $L^1$ and
for any differential ideal in $Q\{y_1,\ldots,y_n\}$ there is a common zero in $L^1$.

Using the same methods we produce the sequence of quasifields $L^n$.
Consider $L^{\omega_1}$, where $\omega_1$ is the first more then countable ordinal.
We shall shaw, that $L=L^{\omega_1}$  is the desired quasifield. Indeed, it
suffices to show that $L$ is differentially closed. Let $\frak u$ be a quasiradical
ideal in $L\{y_1,\ldots,y_n\}$, then from corollary~\ref{GeneratorCor} it follows that
the ideal  $\frak u=\{r_n\}_{n\in\mathbb N}$ is not more than countable generated. Thus
all $r_n$ are determined over some $L^\beta$. Hence, there is a common zero in $L^{\beta+1}$
for all $r_n$. Since the kernel of the homomorphism
$$
L\{y_1,\ldots,y_n\}\to L,
$$
where $y_i$ is mapped to the coordinates of common zero,is a
quasiprime ideal. Then whole ideal $\frak u$ is vanishing on the
mentioned common zero. In other words, we have just checked the
item~(2) of theorem~\ref{DiffClosedFull}.
\end{proof}

It should be noted that in the case of more than countable residue field the existence
of the such sequence of the quasifield  immediately follows from statement~\ref{SubQuasiState}.

\section{$\kappa$-closeness}

Let $Q$ be a quasifield with a residue field $K$. Let fix an
infinite cardinal number $\kappa$. Consider the ring of differential
polynomials $Q\{Y\}$ where $Y$ is an arbitrary set. For any set of
differential polynomials $E$ the set of all its common zeros in
$Q^Y$ will be denoted by $V(E)$, i.~e.
$$
V(E)=\{\,a\in Q^Y\mid \forall g\in E:\: g(a)=0\,\}.
$$
Conversely, for any subset $X$ in $Q^Y$ the set of all polynomials
vanishing on $X$ is denoted by $I(X)$, i.~e.
$$
I(X)=\{\,f\in Q\{Y\}\mid f|_{X}=0\,\}.
$$
It is obvious that for any differential ideal $\frak a$ we have
$\frak{rad}(\frak a)\subseteq I(V(\frak a))$.

The quasifield $Q$ will be said to be a $\kappa$-differentially
closed if for any set $Y$ such that $|Y|\leqslant \kappa$ and any
differential ideal $\frak a$ in $Q\{Y\}$ there is the equality
$\frak{rad}(\frak a)=I(V(\frak a))$.

The following theorem is the natural generalization of
result~\ref{BasicDef} in the case of arbitrary number of variables.

\begin{theorem}\label{DefOfDiffCloseKappa}
The following conditions on quasifield $Q$ are equivalent:
\begin{enumerate}
\item $Q$ is $\kappa$-differentially closed.
\item
For any set $X$ such that $|X|\leqslant\kappa$ and any proper differential
ideal $\frak a$ in $Q\{X\}$ there is a common zero of $\frak a$ in $Q^X$.
\item
Every quasifield generated by not more then $\kappa$ elements over
$Q$ coincides with $Q$.
\end{enumerate}
\end{theorem}
\begin{proof}
It is clear that the conditions~(2) and~(3) are equivalent. Also we see
that~(1) implies~(2). Thus we just only need
to show another implication. Let $\frak a$ be a differential ideal
in $Q\{X\}$ and $f$ is not in $\frak{rad}(\frak a)$. Then there is a
derivative $\theta f$ that is not a nilpotent modulo $\frak a$. Thus
the algebra $Q\{X\}_{\theta f}/\frak a$ is not trivial and therefore
the ideal $[\frak a,z\theta f-1]$ in $Q\{X\cup\{\,z\,\}\}$ is
proper. Since $\kappa$ is infinite then $|X|=|X\cup\{\,z\,\}|$ and
from~(2) it follows that $f$ is not in $I(V(\frak a))$.
\end{proof}

In the same manner as statement~\ref{ClosedHurwitz} and theorem~\ref{DiffClosSingle}
the following result is proved.

\begin{theorem}
A quasifield is $\kappa$-differentially closed iff it is isomorphic
to a Hurwitz series ring over algebraically closed field of
cardinality more than $\kappa$.
\end{theorem}
\begin{proof}
Let show that the mentioned ring of Hurwitz series satisfies the
mentioned properties. Let $K$ be the mentioned algebraically closed
field. Let show that condition~(3) of
theorem~\ref{DefOfDiffCloseKappa} holds.

Indeed, let $L$ be differentially $\kappa$ generated quasifield over
$\Hur K$. The it can be presented in the following form
$$
L=\Hur K\{Y\}/\frak m,
$$
where $\frak m$ is a quasimaximal ideal and $Y$ is of the cardinality less
then or equal to $\kappa$. We shall construct a homomorphism from $L$ to $K$.
Let $\frak n$ be a radical of $\frak m$, then it is clear that it contains
$\Hur K_1\{Y\}$. From theorem~\ref{WeakNullGeom} it follows that
$$
\Hur  K\{Y\}/\frak n= K\{Y\}/\frak n= K.
$$
Let denote the constructed homomorphism $L\to K$ by $\varphi$.  Now
from theorem~\ref{HurwitzUnivers} it follows that there is a unique
differential homomorphism $\Phi$ such that
$$
\xymatrix{
                                       & \Hur K\ar[d]^{\pi}  \\
L\ar@{-->}[ur]^{\Phi} \ar[r]^-{\varphi} & K }
$$
From the uniqueness of the Taylor homomorphism it follows that it
defines over $\Hur K$. Since $L$ is a quasifield then $\Phi$ is an isomorphism.

Let us show the other implication. Let $Q$ be $\kappa$-differentially closed
quasifield with a residue field $K$. It suffices to show that $K$ is of
cardinality more than $\kappa$. Indeed, is it is so, the desired result
immediately follows from statement~\ref{DiffClosSingle}.

Let suppose that contrary holds, in other words, that $K$ is of
cardinality less then or equal to $\kappa$. Then in the polynomial
ring $K[\theta Y]$ there is a maximal ideal $\frak m$ such that its
residue field does not coincides with $K$ (see the proof of
theorem~\ref{WeakNullAlg}). The ideal $\frak m$ can be extended to
the maximal ideal $\frak m'$ in $Q\{Y\}$. Let $\frak q=\pi(\frak
m')$ be a corresponding quasimaximal ideal. Let us show that
quasifield $D=Q\{Y\}/\frak q$ does not coincides with $Q$. For that
we shall compare its residue fields
$$
Q\{Y\}/\frak m'=K[\Theta Y]/\frak m=L.
$$
\end{proof}

We are able to define the notion of a $\kappa$-noetherian
particulary ordered set: a particulary ordered set is called
$\kappa$-noetherian if any ascending chain of its elements of
cardinality more than $\kappa$ is stable. Following the proof of
theorem~\ref{OmegaNoeth} we get the following:

\begin{theorem}
Let $Q$ be an arbitrary quasifield and cardinality of $Y$ is less
then or equal to $\kappa$. Then the set $\QRad Q\{Y\}$ is
$\kappa$-noetherian.
\end{theorem}
\begin{proof}
Let the residue field of $Q$ is denoted by $K$. Let $\frak n$ be a nilradical
of $Q$. We know that the particulary ordered sets
$$
\QRad Q\{Y\}\:\mbox{ and }\:\Rad Q\{Y\}
$$
are isomorphic. The last one is isomorphic to
$$
\Rad Q/\frak n\{Y\}=\Rad K\{Y\}.
$$
We shall show that the set of all ideals in the last ring is $\kappa$-noetherian.
Indeed, Let $Y$ be well-ordered and let $Y_\alpha$  be the subset of all elements
with number less then or equal to $\alpha$. Then every ideal $\frak a$ can be presented
in the following form
$$
\frak a=\bigcup_{\alpha\leqslant\kappa} \frak a_\alpha,
$$
where $\frak a_\alpha=\frak a\cap K\{Y_\alpha\}$.  From the well-known properties of
cardinal numbers it follows that the ideal each ideal $\frak a_\alpha$ is not more then
$\alpha$ generated and thus $\frak a$ is not more than $\kappa$ generated.
\end{proof}

It should be noted that the notion of $\kappa$-differentially closed quasifield
allow us to define the notion of quasivariety in the affine space of any cardinal
dimension and all results about usual quasivarieties can be generalized in to the
case of arbitrary cardinal $\kappa$. But we shall not do it because it is a simple
technical exercise and this question is out of our interests.

\appendix

\section{Appendix}

\subsection{Topology on quasispectrum}

Let $A$ be an arbitrary differential ring, $X=\QSpec A$ be its quasispectrum, and
for each subset $E$ of $A$ we shall define the set
$$
V(E)=\{\,\frak q\in \QSpec A\mid E\subseteq \frak q\,\}.
$$

\begin{statement}\label{QTopology}
In the mentioned notation the following statements holds
\begin{enumerate}
\item If $\frak a$ is a differential ideal generated by the set $E$ then
$$
V(E)=V(\frak a)=V(\frak{rad}(\frak a)).
$$
\item  $V(0)=X$, $V(1)=\emptyset$.
\item Let $\{E_i\}_{i\in I}$ be a family of subsets in $A$. Then
$$
V\left(\bigcup_{i\in I}E_i\right)=\bigcap_{i\in I}V(E_i).
$$
\item $V(\frak a\cap \frak b)=V(\frak a\frak b)=V(\frak a)\cup V(\frak
b)$ for any differential ideals  $\frak a$, $\frak b$ in $A$.
\end{enumerate}
\end{statement}
\begin{proof}
To prove item~(1) it suffices to show the following equality
$$
\frak{rad}(\frak a)=\bigcap_{\frak a\subseteq \frak q} \frak q,
$$
where in the right part all ideals are quasiprime. Items~(2) and~(3)
are obvious. Let show that item~(4) holds. For that it suffices to
show that any quasiprime ideal satisfies the following  property:
from inclusion $\frak a\frak b\subseteq \frak q$ it follows either
$\frak a\subseteq \frak q$, or  $\frak b\subseteq \frak q$. Indeed,
we have the following sequence of inclusions
$$
\frak a\frak b\subseteq \frak q\subseteq \frak r(\frak q).
$$
But $\frak r(\frak q)$ is prime and thus, for example,  $\frak a\subseteq \frak r(\frak q)$.
But $\frak a$ is differential, therefore $\frak a\subseteq\pi(\frak r(\frak q))=\frak q$.
\end{proof}

\subsection{General forms of Nullstellensatz}

\begin{theorem}[Weak form, algebraic variant]\label{WeakNullAlg}
Let $K$ be a field and $\kappa$ is arbitrary cardinal, then
the following are equivalent:
\begin{enumerate}
\item $|K|>\kappa$ or $\kappa$ is finite.
\item Any overfield over $K$ generated by not more than $\kappa$ elements is algebraic over $K$.
\end{enumerate}
\end{theorem}
\begin{proof}
(1)$\Rightarrow$(2). If $\kappa$ is finite then it is a well-known
the Gilbert theorem. Let $\kappa$ be infinite. Let $X$ be the set of
generators of $L$ over $K$. Suppose that contrary holds. Let
$Z\subseteq X$ be a transcendence basis of $L$ over $K$, then
$Y=X\setminus Z$ is algebraic over $K[Z]$. For any element $y\in Y$
there is a polynomial of algebraic dependance over $K[Z]$. Let its
leading coefficient is denoted by $s_y$. Let $S$ be a multiplicative
set generated by all elements $s_y$. Then $L$ is integral over
$S^{-1}K[Z]$. Hence from theorem~\cite[chapter~5, sec.~5,
theor.~5.10]{AM} it follows that $S^{-1}K[Z]$ is a field too. Let
fix one element $z\in Z$, then the set of polynomials $z-\alpha$,
where $\alpha \in K$, is of cardinality more than $\kappa$. The set
$S$ is of cardinality less then or equal to $\kappa$, therefore for
some $\alpha$ the polynomial $z-\alpha$ does not divide the elements
od $S$. Thus $(z-\alpha)$ is nontrivial ideal of $S^{-1}K[Z]$,
contradiction.

(2)$\Rightarrow$(1). Let suppose that $|K|\leqslant\kappa$ and
$\kappa$ is infinite. Then $L=K(x)$, where $x$ is algebraically
independent over $K$, is of cardinality less then or equal to
$\kappa$. Therefore $L$ is not more than $\kappa$ generated over
$K$.
\end{proof}

\begin{theorem}[Full form, algebraic variant]\label{FullNullAlg}
Let $K$ be a field, $Y$ be an arbitrary set of cardinality $\kappa$, then
the following are equivalent:
\begin{enumerate}
\item $|K|>\kappa$ or $\kappa$ is finite.
\item A polynomial ring $K[Y]$ is a Jacobson ring.
\end{enumerate}
\end{theorem}
\begin{proof}
(1)$\Rightarrow$(2). If $\kappa$ is finite then it is a well-known
the Gilbert theorem. Let $\kappa$ be infinite, let $\frak a$  be an
ideal of $K[Y]$, and an element $x$ does not belong to $\frak a$. We
need to show that there exists a maximal ideal $\frak m$ containing
$\frak a$ and not containing $x$. From the choice of the element $x$
it follows that the ring
$$
K[Y]_{x}/\frak a
$$
is nontrivial. Let $\frak m'$ be its maximal ideal. Let its
contraction to $K[Y]$ be denoted by $\frak m$. We need to show that
the ideal $\frak m$ is maximal too. Consider the following sequence
of rings
$$
K\subseteq K[Y]/\frak m\subseteq K[Y]_{x}/\frak m'.
$$
From the previous theorem it follows that the last ring is integral
over $K$. Hence the ring $K[Y]/\frak m$ is integral over $K$. But
$K$ is a field, from theorem~\cite[chapter~5, sec.~2, prop.~5.7]{AM}
it follows that $K[Y]/\frak m$ is a field too.

(2)$\Rightarrow$(1). Item~(2) implies the item~(2) of the previous
theorem, thus the first item holds too.
\end{proof}

\begin{theorem}[Weak form, geometric variant]\label{WeakNullGeom}
Let $K$ be a field, $Y$ be an arbitrary set of cardinality $\kappa$, then
the following are equivalent:
\begin{enumerate}
\item $K$ is algebraically closed and  $|K|>\kappa$.
\item For each nontrivial ideal $\frak a$ of $K[Y]$ there is a common zero in $K^Y$.
\end{enumerate}
\end{theorem}
\begin{proof}
This result immediately follows from algebraic variant because algebraically closed field has
no nontrivial algebraic extensions.
\end{proof}

\begin{theorem}[Full form, geometric variant]\label{FullNullGeom}
Let $K$ be a field, $Y$ be an arbitrary set of cardinality $\kappa$, then
the following are equivalent:
\begin{enumerate}
\item $K$ is algebraically closed and $|K|>\kappa$.
\item For each ideal  $\frak a$ the following holds: $\frak r(\frak  a)=I(V(\frak
a))$, where $V(\frak a)$ is the set of all common zeros of $\frak a$ in $K^Y$, $I(X)$ is the set
of all polynomials vanishing on $X$, and $\frak r (\frak a)$ is a radical.
\end{enumerate}
\end{theorem}
\begin{proof}
From the previous statement it follows that the maximal ideals of $K[Y]$ corresponds to the points of $K^Y$.
Thus our theorem follows from theorem~\ref{FullNullAlg}.
\end{proof}

\end{document}